\documentclass{llncs}

\usepackage{amsmath, graphics, graphicx, float}
\usepackage{verbatim, listings, xspace}
\usepackage{tcolorbox}

\newcommand{\Pd}[3]{\frac{\partial^#3 #1}{\partial #2^#3}}

\newcommand{\opt}[1]{{[\kern-0.225ex[#1]\kern-0.25ex]}}
\definecolor{rowColorA}{HTML}{EDEDED}
\definecolor{rowColorB}{HTML}{FFFFFF}

\def\occa{{\sc occa}\xspace}
\def\OCCA{{\sc OCCA}\xspace}
\def\BOCCA{{\sc {\bf OCCA}}\xspace}

\definecolor{lstBGColor}{HTML}{FFFFFF}
\definecolor{lstBorderColor}{HTML}{606060}
\definecolor{lstTopColor}{HTML}{A1A1A1}

\definecolor{lstKeywordColor}{HTML}{307ED8}
\definecolor{lstCommentColor}{HTML}{D8304C}
\definecolor{lstStringColor}{HTML}{1CBA56}

\definecolor{shadeColor}{HTML}{F5F5F5}

\newcounter{codeBoxCounter}
\newenvironment{codeBox}[1][]{
   \refstepcounter{codeBoxCounter}
   \ifstrempty{#1}
      {\vspace{3mm}\begin{tcolorbox}[breakable,
                                     enhanced,
                                     lines before break,
                                     colback=lstBGColor,
                                     colframe=lstBorderColor,
                                     fonttitle=\sc,
                                     toptitle=1.0mm,
                                     top=-3pt,
                                     bottom=-2pt]}
      {\vspace{3mm}\begin{tcolorbox}[breakable,
                                     enhanced,
                                     lines before break,
                                     colback=lstBGColor,
                                     colframe=lstBorderColor,
                                     fonttitle=\sc,
                                     top=-3pt,
                                     bottom=-2pt,
                                     toptitle=1.0mm,
                                     bottomtitle=1.0mm,
                                     title=\thecodeBoxCounter.\ #1]}
}
{\end{tcolorbox}}

\newcommand{\sCodeBox}[1][]{
   \ifstrempty{#1}
   {\begin{codeBox}}
   {\begin{codeBox}[#1]}
}
\def\eCodeBox{\end{codeBox}}

\def\lstSidePadding{2mm}
\def\lstSideMargin{3mm}
\def\lstNumberPadding{4mm}

\lstdefinelanguage{occa}[ISO]{C++}{
   keywords={
      occa, occaMemory,
      global, local, barrier,
      occaInnerDim0,occaInnerDim1,occaInnerDim2,
      occaOuterDim0,occaOuterDim1,occaOuterDim2,
      occaInnerFor,occaInnerFor0,occaInnerFor1,occaInnerFor2,
      occaOuterFor0,occaOuterFor1,occaOuterFor2,
      occaGlobalFor0,occaGlobalFor1,occaGlobalFor2,
      occaGlobalId0,occaGlobalId1,occaGlobalId2,
      occaGlobalDim0,occaGlobalDim1,occaGlobalDim2,
      occaInnerId0,occaInnerId1,occaInnerId2,
      occaOuterId0,occaOuterId1,occaOuterId2,
      occaLocalMemFence,occaGlobalMemFence,
      occaShared,occaPointer,occaConstant,occaVariable,occaRestrict,occaVolatile,occaConst,occaAligned,
      occaKernelInfoArg,occaFunctionInfArg,occaFunctionInfArg,occaFunctionInfArg,occaFunctionInfArg,occaFunctionInfo,
      occaKernel,occaFunction,occaDeviceFunction,
      occaFunctionShared,occaPrivateArray,occaPrivate,
      occaUnroll,occaInnerReturn}
}

\def\codedef#1,#2,[#3-#4],#5\relax{
    \lstinputlisting[title=#1,
                     caption=#1,
                     firstnumber=#3,
                     linerange={#3-#4},
                     language=#5,
                     linewidth={\linewidth\dimexpr\linewidth\relax},
                     numbers=left,
                     numberfirstline=false,
                     numbersep=\lstNumberPadding,
                     backgroundcolor=\color{lstBGColor},
                     showtabs=false,
                     frame=single,
                     tabsize=2,
                     captionpos=t,
                     basicstyle=\small\ttfamily,
                     keywordstyle=\color{lstKeywordColor},
                     commentstyle=\color{lstCommentColor},
                     stringstyle=\color{lstStringColor},
                     breaklines=true,
                     breakatwhitespace=false,
                     showstringspaces=false,
                     framexleftmargin=\lstSidePadding,
                     framexrightmargin=\lstSidePadding,
                     xleftmargin=\lstSideMargin,
                     xrightmargin=\lstSideMargin,
                     lineskip=0.3mm,
                     escapeinside={{*@}{@*}},
                     morekeywords={*,...}]{#2}
}

\newcommand\setCode[1]{\setCodedef#1\relax}
\def\setCodedef#1,#2,#3\relax{
  \lstset{
    language=#1,
    linewidth={\dimexpr\linewidth\relax},
    numbers=#2,
    numberfirstline=false,
    numbersep=\lstNumberPadding,
    backgroundcolor=\color{lstBGColor},
    showtabs=false,
    frame=single,
    tabsize=2,
    captionpos=b,
    basicstyle=\ttfamily#3,
    keywordstyle=\color{lstKeywordColor},
    commentstyle=\color{lstCommentColor},
    stringstyle=\color{lstStringColor},
    breaklines=true,
    breakatwhitespace=false,
    showstringspaces=false,
    framexleftmargin=\lstSidePadding,
    framexrightmargin=\lstSidePadding,
    xleftmargin=\lstSideMargin,
    xrightmargin=\lstSideMargin,
    lineskip=0.3mm,
    escapeinside={{*@}{@*}},
    morekeywords={*,...},
  }
}

\def\setCodeWithTitledef#1,#2,#3,#4\relax{
  \lstset{
    title=#1,
    language=#2,
    linewidth={\dimexpr\linewidth\relax},
    numbers=#3,
    numberfirstline=false,
    numbersep=\lstNumberPadding,
    backgroundcolor=\color{lstBGColor},
    showtabs=false,
    frame=single,
    tabsize=2,
    captionpos=b,
    basicstyle=\ttfamily#4,
    keywordstyle=\color{lstKeywordColor},
    commentstyle=\color{lstCommentColor},
    stringstyle=\color{lstStringColor},
    breaklines=true,
    breakatwhitespace=false,
    showstringspaces=false,
    framexleftmargin=\lstSidePadding,
    framexrightmargin=\lstSidePadding,
    lineskip=0.3mm,
    escapeinside={{*@}{@*}},
    morekeywords={*,...},
  }
}

\title{High-Order Finite-differences on multi-threaded architectures using \OCCA \thanks{Partial funding from Royal Dutch Shell, ONR award N00014-13-1-0873, and sub-contract to the Argonne CESAR Exascale Co-design Center under award number ANL 1F-32301. Shell release of technical information: RTI 2094414 \& 233820614.}}
\titlerunning{\OCCA}

\author{David S. Medina\inst{1} \and Amik St-Cyr\inst{2} \and Timothy Warburton\inst{1}}
\institute{Rice University, Computational and Applied Mathematics, \email{\{dsm5,tim.warburton\}@rice.edu}
\and Royal Dutch Shell, Seismic applications team, \email{amik.st-cyr@shell.com}}

\begin{document}

\maketitle

% Already too long
\begin{abstract}
High-order finite-difference methods are commonly used in wave propagators for industrial subsurface imaging algorithms. Computational
aspects of the reduced linear elastic vertical transversely isotropic propagator
are considered. Thread parallel algorithms suitable for implementing this propagator on multi-core and many-core
processing devices are introduced. Portability is addressed through the use of the \OCCA
runtime programming interface. Finally, performance results are shown for various architectures
on a representative synthetic test case.
\end{abstract}

\section{Introduction}
%
% TIM: please write a short intro (or go straight to the problem description?)... Also some trimming of the abstract?
%

High-order finite-differences are used in seismic imaging and many other industrial applications primarily
because of their computational efficiency.
High-order wave propagators lie at the heart of numerous seismic imaging applications,
 such as full waveform inversion and reverse time migration (RTM).
We study multi-threaded performance on various current and emerging computing architectures of
propagators for vertical transversely isotropic media (VTI).

The VTI propagator introduced in \cite{du08} is given by
\begin{align}
\Pd{p}{t}{2} &= \nu^2_x
\left[
\Pd{p}{x}{2}+
\Pd{p}{y}{2}\right] +
\nu^2_z
\Pd{q}{z}{2} + s(t)\delta({\bf x} - {\bf x_i})\label{eq:p}, \\
\Pd{q}{t}{2} &= \nu^2_n
\left[
\Pd{p}{x}{2}+
\Pd{p}{y}{2}\right] +
\nu^2_z
\Pd{q}{z}{2}\label{eq:q}.
\end{align}
In the preceding equations, $p$ is an approximation for the $P$-wave
while $q$ is and auxiliary wavefield variable. $\epsilon$ and $\delta$
are the anisotropic parameters. The vertical $P$-wave velocity is
represented with $\nu_z$ and its horizontal component is
$\nu_x =\nu_z \sqrt{1+2\epsilon}$ while the normal move-out velocity
is $\nu_n=\nu_z \sqrt{1+2\delta}$. For this approximation to be
relevant $\epsilon-\delta \le 0$ is necessary. The forcing considered
in our benchmark is the Ricker wavelet $s=(1-2\pi^2 f^2 t^2) e^{-\pi^2 f^2 t^2}$ with $f=15Hz$.

We consider a centered finite-difference discretization in time and space in second order form on infinite domains.
For ${\bf u}({\bf x},t)=(p,q)^T$ and ${\bf F}({\bf u},{\bf x},t) $ set as the right and side of (\ref{eq:p})-(\ref{eq:q})
the centered in time approximation reads
\begin{align}
{\bf u}^{n+1} - 2{\bf u}^{n} + {\bf u}^{n-1} \approx \Delta t^2 {\bf F}({\bf u}^{n})
\end{align}
where ${\bf u}^k \equiv {\bf u}({\bf x},t^n)$ with $t^n=n \Delta t$.
High-order finite-difference stencils are of practical importance for the
efficient numerical solutions of wave propagation problems \cite{alford74,vishnevsky14}.
Indeed, for a similar number of points composing the computational grid, the
number of points required to resolve the shortest wavelength (as defined by Nyquist)
decreases and gets close to the spectral or pseudo-spectral limit of two
points per wavelength \cite{fornberg87}.
Most propagators used in seismic applications use two different flavors of
high-order finite-differences. The earth subsurface is geologically horizontally
layered. Since depth, represented by the $z$ coordinate, will experience the
most changes in the rock properties, while in the $x-y$ planes the properties
will remain constant within a layer. Therefore, a common strategy is to
have a symmetric stencil in the $x-y$ direction, while handling a variable
spacing in $z$. The weights and spacings can be optimized
to handle a variety of physical and numerical properties
\cite{holberg87,fornberg98}.
For simplicity, we suppose a domain $\Omega=[0,Lx]\times[0,Ly]\times[0,Lz]$
where $\Delta x = \Delta y = h$ and $\Delta z_k$ result from the discretization in space
using $N_{d=\{x,y,z\}}$ points in each direction respectively. The mesh size in the $z$ direction
varies per grid point belonging to a different $x-y$ plane. Adopting the convention $p(x_i,y_j,z_k)= p_{i,j,k}$,
the differentiation stencil in the $x-y$ plane is
\begin{equation}
h^2(\frac{\partial^2 }{\partial x^2} + \frac{\partial^2 }{\partial y^2})p_{i,j,k}
\approx
w^{xy}_0 p_{i,j,k}+\sum^{R_{xy}}_{l=1} w^{xy}_l\left(p_{i+l,j,k}+p_{i-l,j,k}+p_{i,j+l,k}+p_{i,j-l,k}\right)
\end{equation}
where the $w^{xy}_l$ are the $R_{xy}+1$ weights for approximating the two dimensional Laplacian.
The differentiation is a bit simpler in the $z$ direction:

\begin{align}
\frac{\partial^2 }{\partial z^2}q_{i,j,k}&\approx \sum^{R_{z}}_{l=-R_{z}+} w^{z}_{k,l} q_{i,j,k+l}.
\end{align}
Again, the $w^{z}_{k,l}$ are the weights for approximated the second derivative. However, for each
position $z_k$ where the value of the derivative is sought, $2R_z+1$ weights are needed instead of $R_z+1$
as in the symmetric case due to the asymmetry in the $z$ direction. The grid size $\Delta z_k$ is absorbed into the $w^{z}_{k,l}$ weights in practice
and therefore are not appearing above. The domain $\Omega$ is embedded into a larger domain where a
damping formula is applied as in \cite{cerjan85}. Outside the damping region, the solution is assumed to
be zero for $R_{xy}$ points in the $x$ and $y$ directions and $R_z$ points in $z$.

In the following sections, we describe the reduced elastic VTI model for isotropic media together with a typical finite-difference discretization employed in industry.

\section{Computational efficiency of high-order finite differences}
\label{sec:ci}
The peak parallel floating point operations per second (flops) available on modern CPUs and GPUs
have followed the trend set by Moores law. Unfortunately, the available memory bandwidth
lagged this trend. This gap in bandwidth currently favors algorithms generating
lots of flops per byte of data moved \cite{mccalpin07}.
For VTI, using this type of stencil, a pessimistic computational intensity is
\begin{align}
CI \approx (1/4)(5R_{xy}+4R_z)/(4R_{xy}+2R_z) \approx 0.4 \;\;\mbox{flops/byte}
\end{align}
where most of the loads are assumed to be not in cache. An idealized version is to consider
the least loads as possible (assumes most of the data in cache). This is done by assuming three
single precision loads per point for the model properties  ($\nu^2_x$, $\nu^2_n$ and $\nu^2_z$) as well
as the two pairs of loads and stores for ${\bf u}^{n}$ and, respectively, ${\bf u}^{n+1}$.
\begin{align}
CI \approx (1/28)(5R_{xy}+4R_z) \approx 0.3 (R_{xy}+R_z) \;\;\mbox{flops/byte}
\end{align}
Therefore increasing the order of the stencil augments the intensity since the low order
case is close to the pessimistic estimate. In practice, better approximations can be obtained
\cite{williams09}. Performing those measurements automatically using hardware counters
is still in development. % [cite]
Moreover in \cite{fornberg87}, the effectiveness of finite-differences for wave propagation
problems is shown to increase with order. Indeed for a fixed number of Fourier modes
``M'', $\tilde{N}_d$ points are required to guarantee their resolution according to Nyquist.
The relation is approximated with $\tilde{N}_d = c_p M^{1+(2R)^{-1}}$ and therefore doubling
the polynomial order for a fixed number of modes, leads to $M^{\frac{1}{4}}$ times more points in each
direction. Since the method is explicit in time, the total increase in computational cost is
$M$ in 3D.

\section{\OCCA: portable multi-threading}

\OCCA, a recently developed C++ library for handling multi-threading is employed. It
uses run-time compilation and macro expansions which results in a novel and simple single
kernel language that expands to multiple threading languages. \OCCA currently supports device kernel expansions
for the OpenMP, OpenCL, pThreads, Intel COI and CUDA languages. Performance characteristics are given for our implementations built on top of the \OCCA API. Using the unified \OCCA programming approach allows customized kernels optimized for CPU and GPU architectures with a single ``host'' code.

\begin{figure}[hb!]
\begin{center}
\includegraphics[width=8cm]{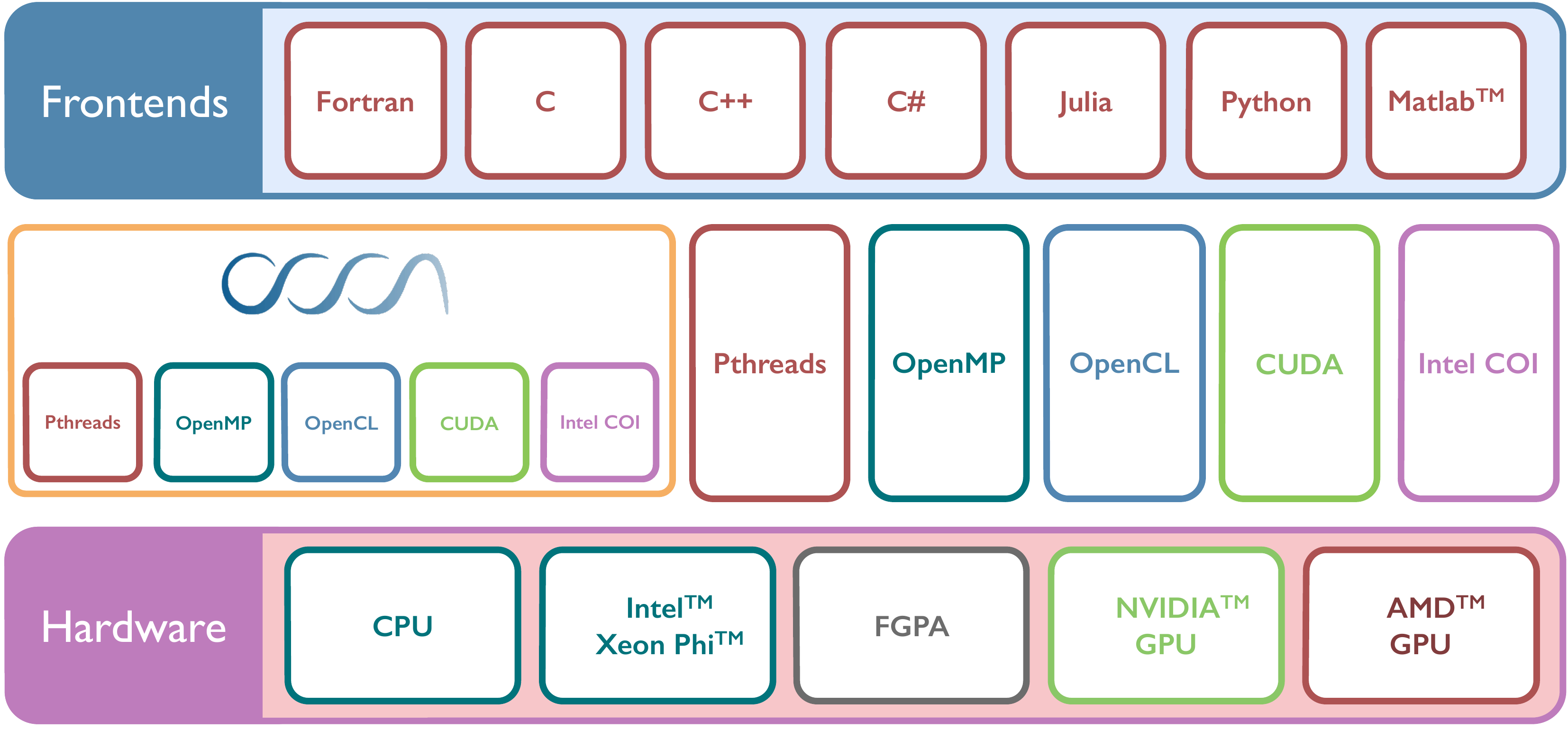}
\end{center}
\caption{\OCCA wraps different language APIs and is non-conflicting with external libraries in either platform}
\label{fig:flow}
\end{figure}

The \occa library is an API providing a kernel language and an abstraction layer to back-ends
APIs such as OpenMP, OpenCL and CUDA see \cite{openmp08,stone10,nickolls08} amongst others.
Although \occa unites different threading platform back-ends, the main contributions is the abstraction
of the kernel language. Using a macro-based approach, an \occa kernel can be expanded at runtime
to suit OpenMP with dynamic pragma insertions or a device kernel using either OpenCL or CUDA and is
taylored to streamline the incorporation of future threading languages with ease.
This is similar to source to source front-ends to compilers (cite HMPP and PGI) however,
the translation is at runtime and performed mostly by the c-preprocessor. With such an approach,
the control over the optimization is completely left to the user and more information is available
for the compiler.
\newline\newline\noindent {\bf \BOCCA host API}:
Aside from language-based libraries from OpenMP, OpenCL or CUDA, the
\occa host API is a stand-alone library.This independence allows \occa to be combined with other
libraries without conflict, as shown in \ref{fig:flow}.
The three key components that influenced the \occa host API development: the platform device,
device memory and device kernels. Presenting the entire \OCCA API is not feasible in this paper.
For the complete details see \cite{medina14} and the git repository for the latest
developments\footnote{\url{http://github.com/tcew/OCCA}}. We try here to expose the minimal
knowledge required to write the VTI kernel.
\newline\newline\noindent {\bf \BOCCA device}:
An \occa device acts as a layer of abstraction between the \occa API and the API from supported
languages. Due to the just-in-time code generation, the platform target can be chosen at run-time.
Enabled platforms are managed at compile-time in the case of unsupported platforms on the
compiled architecture.
An \occa device generates a self-contained context and command queue \opt{stream} from a chosen
device, being a socketed processor, GPU or other OpenCL supported devices such as a Xeon Phi
or an FPGA. Asynchronous computations with multiple contexts can be achieved using multiple
\occa devices. A device object can allocate memory and compile kernels for the physical device.
\newline\newline\noindent {\bf \BOCCA memory}:
The \occa memory class abstracts the different device memory handles and provide some useful
information such as device array sizes. Although memory handling in \occa facilitates host-device
communication, the management of reading and writing between host and device is left to the programmer
for performance reasons. The dedicated device memory class allows the
\occa kernel class to manage communications between distinct memory types.
\newline\newline\noindent {\bf \BOCCA kernel}:
The \occa kernel class unites device function handles with a single interface, whether for a function
pointer (OpenMP), cl\_kernel (OpenCL), or cuFunction (CUDA). When using the OpenCL and CUDA
kernel handles, passing the arguments through their respective API is simple. These implicit
work-group \opt{block} and work-item \opt{thread} sizes are passed by an argument in CPU-modes such as OpenMP.
\newline\newline\noindent {\bf \BOCCA kernel language}:
 GPU computing involves many threads and the thread-space is
logically decomposed into thread-blocks. Thread blocks are queued for execution onto the available multiprocessors.
In general a GPU chip has more than a single multiprocessor and the choices for number of blocks and threads per blocks
are dependent on the algorithm, resources available and the developer. Taking as example the CUDA API, a general loop
can be written as in {\bf\small Listing \ref{lst:cuLoopExpansion}}.
\setCode{occa,none,\scriptsize}
\begin{lstlisting}[caption={
  The expansion of the implicit for-loops found in CUDA kernels is displayed.
  Loop grouping (1) expands multi-dimensional work-groups \opt{blocks} and loop grouping (2)
  expands multi-dimensional work-items \opt{threads}.
},label={lst:cuLoopExpansion}]
for(int bZ = 0; bZ < gridDim.z; ++bZ){        // (1)
  for(int bY = 0; bY < gridDim.y; ++bY){
    for(int bX = 0; bX < gridDim.x; ++bX){
      // Shared memory is initialized here
      for(int tZ = 0; tZ < blockDim.z; ++tZ){  // (2)
        for(int tY = 0; tY < blockDim.y; ++tY){
          for(int tX = 0; tX < blockDim.x; ++tX){
             // Work here,  initialize register memory
}}}}}}
\end{lstlisting}
The first three for loops are going through all the blocks while the three innermost loops
iterate on the individual threads per block.This is an implicit loop, one would never write such
a thing in practice as it would destroy the thread concurrency and all achievable performance.
\OCCA replaces the loops with the correct constructs and gives the developer access to the
global indexing using {\tt\small occaGlobalId\{0,1,2\}} which is computed with {\tt\small blockIdx}
and {\tt\small threadIdx}.
The \occa equivalent general loop is contained in listing \ref{lst:occaLoopExpansion}.
\vspace{4mm}
\setCode{occa,none,\scriptsize}
\begin{lstlisting}[caption={
  The \occa programming model mirrors GPU programming, where group loopings (1) and
  (2) refer to work-groups \opt{blocks} and work-items \opt{threads} respectively.
},label={lst:occaLoopExpansion}]
occaOuterFor2{       // Loop grouping (1)
  occaOuterFor1{
    occaOuterFor0{
      // Shared memory defined here
      occaInnerFor2{ // Loop grouping (2)
        occaInnerFor1{
          occaInnerFor0{
            // work here
}}}}}}
\end{lstlisting}

\noindent The use of shared memory is still available in \OCCA since it is essential for many GPU optimized codes.
Shared memory still acts as a scratchpad cache for GPU architectures but can be seen as a prefetch buffer for CPU-modes in \OCCA.

\begin{figure}[hb!]
\begin{center}
\includegraphics[width=4.5cm]{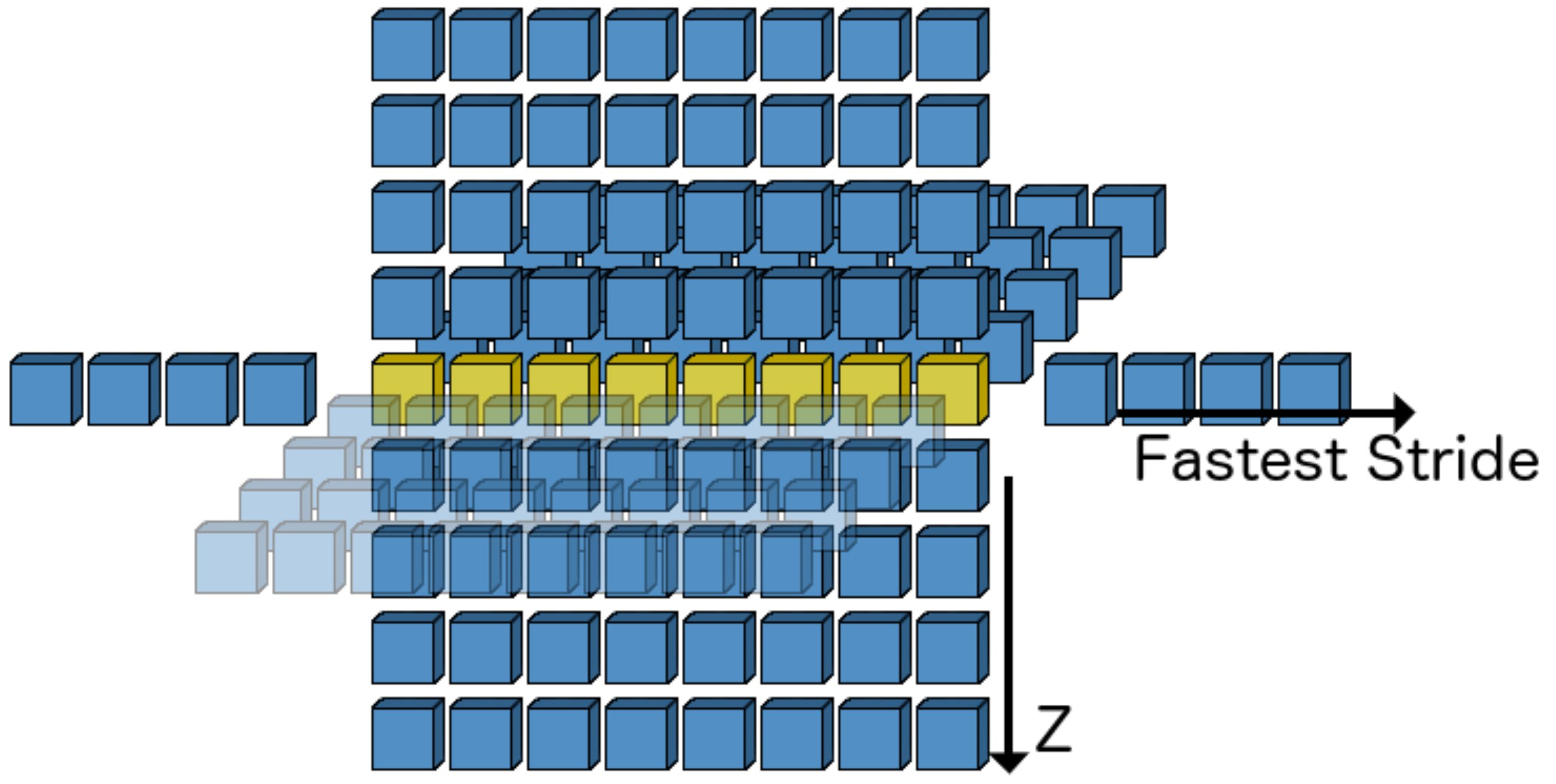}
\hspace{10mm}
\includegraphics[width=4.5cm]{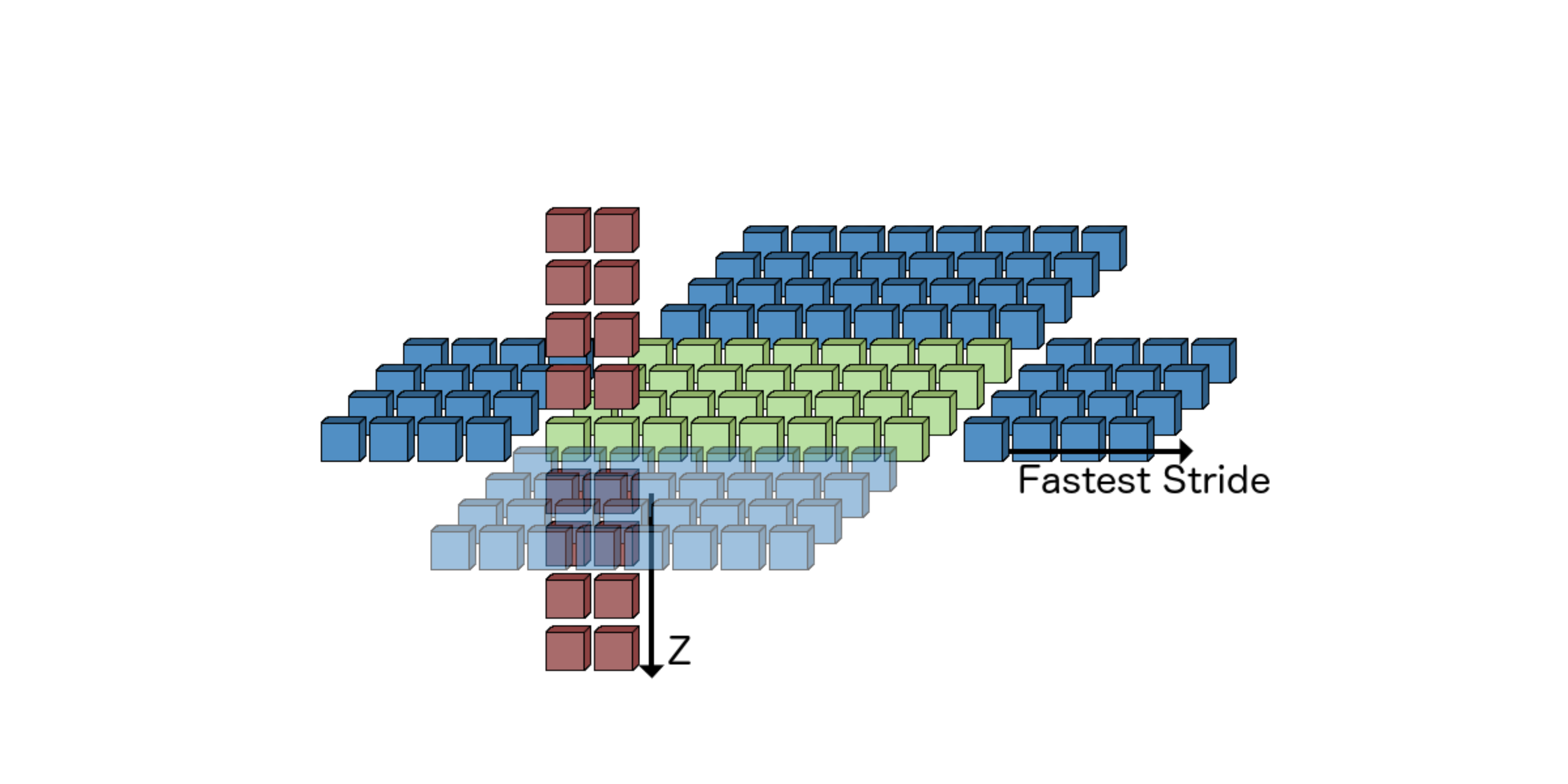}
\end{center}
\caption{The left panel represents a 3D finite-difference stencil vectorized with AVX. The fast stride is in the x direction and $8$ single precision
stencil evaluation are performed simultaneously. The right panel represents thread block with the large 2D subdomain the information loaded
into shared (fast) memory. The register rolling in the z-direction is shown and for a "two-elements" kernel, each thread handles two columns.}
\label{fig:cpu_gpu}
\end{figure}
The \OCCA:OpenMP code performs the VTI steps uses the classic technique of cache blocking as seen in code listing \ref{lst:vtiPseudocode}.
The best performing kernel had 2D cache blocking with the $Z$-block first followed by the $Y$-block.
The innermost loop would be $x$ (stride-$1$) then $z$ and finally $y$. The z-blocks were handled in
\occa with {\tt\small occaOuterFor2}, the $y$ one with {\tt\small occaOuterFor1} and so on. The
vectorization was handled directly by the Intel compiler by placing a pragma {\tt\small \#pragma ivdep}
in the {\tt\small occaInnerFor0} ($x$) and making sure the data was correctly padded.
The size of the blocking in $z$-$y$ was determined as $(28,20)$ by
running the code over a set of grids and possible block ranges and comparing throughput times see section \ref{sec:perf}.
The OpenMP first-touch policy was critical in obtaining performance across dual sockets as well as the
correct thread affinity. Finally, to make sure the compiler was optimizing as depicted in Fig. \ref{fig:cpu_gpu},
a hand written kernel with explicit register blocking was written: $5\%$ increase in performance
was observed.

\begin{lstlisting}[caption={For each time-step, the 2D blocks at the top of the structured grid sweep in the $z$ direction and update all points in the current $z$ plane.},label={lst:vtiPseudocode}]
Partition the top plane of the grid into *@$B_x\times B_y$@* blocks of size *@$w\times h$@*

For time-step *@$n = 0, 1, \ldots,$@* time-Steps
    For each block *@$(b_i,b_j)$@*                                   (1)
        For *@$n = 0, 1, \ldots, N_z$@*
            For each point *@$(i,j,k)$@* such that
                (*@$b_i \leq i < b_i + w$@*) and (*@$b_j \leq j < b_j + h$@*)   *@$\hphantom{a}$@*         (2)
                Update *@$p^{n+1}(i,j,k)$@* and *@$q^{n+1}(i,j,k)$@*
            End For // Point Update
        End For     // Traversing depth
    End For         // Iterating over blocks
End For             // Computing a time-step update
\end{lstlisting}

A single implementation encompasses \OCCA:OpenCL and \OCCA:CUDA follows directly the work of \cite{paulius}. As depicted in the right
panel of Fig. \ref{fig:cpu_gpu}, for a given thread block, the 2D $x-y$ stencil executes into fast shared memory
while the $z$ direction is handled by register rolling. If each thread handles one such column per thread block then
this is a {\em one-element} approach while a {\em two-elements} approach consists of having two such columns
per thread. Care was taken to align the data to enable coalescing loads to shared memory.

\section{Performance}
\label{sec:perf}
The VTI kernel is integrated in time for a thousand time steps. A metric of performance used in seismic is the
throughput: number of sweeps through the entire grid block per second. The precision is set at $R_{xy}=12$
and $R_z=8$ and yields approximatively $92$ flops per point. The CI optimistic model derived in section \ref{sec:ci}
yields a factor of $3.3$.

Results on a dual socket node with E$5$-$2670$ are reported in table \ref{table:cpu}. The dual socket node
is capable of $666$ single precision GFlops while the bandwidth is $102.4$ GB/s.
The optimistic CI predicts a maximal peak of 47\%. The results show the fastest \occa
kernel achieving 21\% and good scalability as compared to the native OpenMP code (without \occa).
The difference stems from the added knowledge at compile time for \occa, where all loop-bounds are known at compile time.

% Show perf OpenMP vs OCCA OpenMP
\begin{table}[h!]
\begin{center}
\begin{tabular}{lcccccc|c}
  \hline
  Project   & Distribution & 1 Thread & 2 Threads   & 4 Threads  & 8 Threads  & 16 Threads  & \% Peak \\ \hline
  Native      & Compact             & 92       & 183 (98\%)  & 360 (96\%) & 668 (89\%) & 1226 (82\%) & 17\\
  Native      & Scatter             & 92       & 183 (98\%)  & 356 (95\%) & 686 (92\%) & 1191 (80\%) & 16\\
  \OCCA & Compact             & 115      & 229 (99\%)  & 448 (97\%) & 820 (89\%) & 1548 (84\%) & 21\\
  \OCCA & Scatter             & 115      & 230 (100\%) & 454 (98\%) & 884 (96\%) & 1411 (76\%) & 19\\
  \hline
\end{tabular}
\end{center}
\caption{Multithreading scaling with OpenMP using alternative thread distributions on different
number of cores (using two Xeon E5-2640 Processors)}\label{table:cpu}
\end{table}

Table \ref{tab:gpuComp} contains performance on GPU architectures that were based on optimized CUDA code and translated to \occa.
We note that performance seen in table \ref{tab:gpuComp} was on par with native code due to optimizations that can be done with run-time compilation including manual unrolling and manual bounds on OpenMP-loops.

\begin{table}[!h]
\begin{center}
\begin{tabular}{llll}
  \hline
  Project \hspace{2cm} & Kernel Language \hspace{1cm} & K10 (1-chip) \hspace{1cm} & K20x  \\ \hline
  Native          & CUDA            & 1068              & 1440 \\
  Native (2)      & CUDA            & 1296              & 2123 \\
  \OCCA           & \OCCA:CUDA      & 1241              & 1934 \\
  \OCCA (2)       & OCCA:CUDA       & 1579              & 2431 \\
  \OCCA           & \OCCA:OpenCL    & 1303              & 1954 \\
  \OCCA (2)       & OCCA:OpenCL     & 1505              & 2525 \\
  \hline
\end{tabular}
\end{center}
\caption{
  Performance comparisons on the VTI update kernels tailored for GPU architectures.
  Update kernels use 1-point updates per work-item/thread or are labeled with (2) to represent 2-point update kernels.
  One K10 chip runs at 745 MHz and contains 1536 floating point units with 160 GB/s bandwidth.
  By comparison, the K20x runs at 732 MHz and contains 2496 floating point units with 250 GB/s bandwidth.\vspace{-4mm}
}
\label{tab:gpuComp}
\end{table}
Table \ref{tab:cpuGpuComp} contains results from two optimized kernels, a CPU-tailored code and a GPU-tailored code, run on OpenMP, OpenCL and CUDA to note performance portability.
Although it was expected that optimal CPU-tailored algorithms would not give optimal performance for GPU architectures, we see 40-50\% of optimal performance by just running the \OCCA kernels in GPU-modes.
The GPU-tailored algorithm ran on CPU-modes ended running on 20\% performance compared with optimal CPU code, mainly due to the lack of direct control over shared memory as seen on GPU architectures.

\begin{table}[!h]
\begin{center}
\begin{tabular}{lll}
  \hline
  \hspace{4cm}        & CPU-tailored Kernel \hspace{1cm} & GPU-tailored Kernel \\ \hline
  OpenMP              & 1548                & 364 (23\%)          \\
  CUDA   (1 K10 core) & 515 (41\%)          & 1241                \\
  OpenCL (1 K10 core) & 665 (51\%)          & 1302                \\
  \hline
\end{tabular}
\end{center}
\caption{
  Performance comparisons between combinations of OpenMP, CUDA and OpenCL running on the CPU and GPU tailored kernels.
}
\label{tab:cpuGpuComp}
\end{table}
\vspace{-9mm}

\section{Conclusion \& future work}
We have studied a vertical transverse isotropic propagator discretized with centered finite-differences
in time and space. Finite-differences are extensively used in seismic modeling. We have justified the
advantage of using high-order stencils both in terms of computational efficiency and points needed
per wavelength. To enable the study on various compute architectures, a multi-threaded gateway API to many
multi-threading APIs was employed: \OCCA. The performance results obtained with the library are generally
faster than with the codes written using the best API for the hardware, thanks to the just-in-time compilation.
For now, it seems a single \OCCA kernel solution performing well for two types of architecture is impossible.
The main factor preventing portable optimization is due to the lack of direct control over cache on CPU architectures which
can be done on GPU architectures through shared memory.
This level of control is currently only available for GPGPUs and unavailable for traditional CPUs. Having such control
on the next generation of CPUs would most certainly re-open possibilities of a single code performing
efficiently on both architectures.

We are currently working on implementing features that will make \OCCA more accessible to programmers.
The current \OCCA project is focusing on the \OCCA Kernel Language (OKL) derived from compiler-tools to allow for a more native kernel language as opposed to the current macro-based language.
A full parser is used for adding language features such as loop-interchange flexibility, automatic OpenCL-CUDA translation to OKL, kernel splitting from multiple outer-loops with the goal to embedded the OKL kernels in regular code.

\bibliography{references}
\bibliographystyle{plain}

% Appendix
% Q: add the code listings?

\end{document}